\documentclass[12pt,twoside]{extarticle}
\usepackage{mathtext}
\usepackage{titlesec}
\usepackage[cp1251]{inputenc}
\usepackage[T2A]{fontenc}
\usepackage[russian,english]{babel}
\usepackage{fancyhdr}
\usepackage{amsthm}
\usepackage{dsfont}
\usepackage{indentfirst}
\usepackage{mathrsfs}
\usepackage{array}
\usepackage{amssymb}
\usepackage{amsmath}
\usepackage{eucal}
\usepackage{dubloper}
\usepackage{bm}
\usepackage{epsfig}
\usepackage{wrapfig}

\titleformat{\section}[block]{\large\bfseries\filcenter}{\large\bfseries\thesection. }{0pt}{}
\titleformat{\subsection}[block]{\bfseries\filcenter}{\bfseries\thesubsection. }{0pt}{}
\titleformat{\subsubsection}[block]{\bfseries\filcenter}{\bfseries\thesubsubsection. }{0pt}{}

\newcounter{mytheorem}[section]

\edef\lim{\lim\limits}
\renewcommand{\frac}[2]{\dfrac{#1}{#2}}
\numberwithin{equation}{section}
\numberwithin{figure}{section}
\begin{document}

\setlength{\abovedisplayskip}{2mm} 
\setlength{\belowdisplayskip}{2mm} 
\parindent=8,8mm
\renewcommand{\headrulewidth}{0pt}
\edef\sum{\sum\limits}

\thispagestyle{empty}

\begin{center}

\

\

\

\textbf{\Large Finite groups with $\mathfrak{F}$-subnormal }

\textbf{\Large normalizers of Sylow subgroups}

\bigskip

\textbf{\bf A.\,F.\,Vasil'ev, T.\,I.\,Vasil'eva, A.\,G.\,Melchenko}

\bigskip

\end{center}

\bigskip

\begin{center}
{\bf Abstract}
\end{center}

\medskip

\noindent\begin{tabular}{m{8mm}m{138mm}}

&{\small~~~
Let $\pi$ be a set of primes and $\mathfrak{F}$ be a  formation. In this article a properties of the class ${\rm w}^{*}_{\pi}\mathfrak{F}$ of all groups $G$,  such that
$\pi(G)\subseteq \pi(\mathfrak{F})$  and the normalizers of all Sylow $p$-subgroups of  $G$ are $\mathfrak{F}$-subnormal in $G$
 for every $p\in\pi\cap\pi(G)$ are investigated.
It is established that ${\rm w}^{*}_{\pi}\mathfrak{F}$ is a formation.
Some hereditary saturated formations $\mathfrak{F}$ for which ${\rm w}^{*}_{\pi}\mathfrak{F}=\mathfrak{F}$ are founded.

\smallskip

}

\smallskip

MSC2010\  20D20, 20D35, 20F16


\end{tabular}

\medskip

{\parindent=0mm
\textbf{\bf {\large Introduction}}
}
\bigskip

We consider only finite groups. It is well known what role is played the properties of normalizers of the primary subgroups (local subgroups) in classification of finite simple non-abelian groups. In recent years, local subgroups are actively used in the study of non-simple, in particular, soluble groups.
 In 1986 it was established [1] that a group is nilpotent if the normalizers of its Sylow subgroups (briefly, Sylow normalizers) are nilpotent.
 Groups with supersoluble Sylow normalizers were studied in [2-4]. A series of papers [5-9] is dedicated to the study of groups whose all the Sylow normalizers belong to a saturated formation $\mathfrak{F}$.

In this paper, we are interested in the following question.
How do the properties of embedding of Sylow normalizers into a group influence on the structure of the whole group?

We note the following results. Group $ G $
is nilpotent if and only if its any Sylow
normalizer coincide with $G$. By the well-known Glauberman's theorem [10], if all Sylow subgroups of a group are self-normalizing, then the group is a $p$-group for some prime $p$.

  Let $H$ be a subgroup of a group $G$. Consider a chain of subgroups

\medskip

\rightline {$H=H_0 \leq H_1 \leq \dots \leq \ H_{n-1} \leq H_{n} = G.\ \ \ \ \ \ \ \ \ \ \ \ \ \ \ \ \ \ \ \ \ \ \ \ \ \ \ \ \ (1)$}
\medskip

According to [11], $H$ is called {\it $\mathbb{P}$-subnormal} in $G$ if either $H = G$ or there exists a chain (1) such that
$|H_{i}:H_{i-1}|$ is a prime for any $i = 1, \dots, n$;

According to [12], $H$ is called {\it ${\rm K}$-$\mathbb{P}$-subnormal} in $G$ if there exists a chain (1) such that either $H_{i-1}\unlhd H_{i}$, or $|H_{i}:H_{i-1}|$ is a prime for any $i = 1, \ldots, n$.

In [13] V.S. Monakhov and V.N. Kniahina established that a group $G$ is supersoluble if and only if all its Sylow normalizers  are $\mathbb{P}$-subnormal in $G$.

A subgroup  $H$ is called {\it submodular} in  $G$ [14], 
if there exists a chain of subgroups
(1) such that $H_{i-1}$ is a modular subgroup in $H_{i}$ for $i=1,\ldots,s$.
Here the {\it modular} subgroup in $G$ is a modular element in the lattice of all subgroups of $G$ [15].
The class $s\mathfrak{U}$ of all strongly supersoluble groups was studied in [16]
($s\mathfrak{U}$ is the class of supersoluble groups,
in which all Sylow subgroups are submodular). By [17, Theorem 3.2], 
if the normalizers of all Sylow subgroups of a group $G$ are submodular, then $G\in s\mathfrak{U}$.

The concept of subnormality was generalized by T.O.~Hawkes [18], L.A.~Shemetkov [19] as follows.

Let $\mathfrak{F}$ be a non-empty formation. A subgroup $H$ is called 
{\it $\mathfrak{F}$-subnormal in $G$}
(which is denoted by $H$
$\mathfrak{F}$-${\rm sn}$ $G$),
if either $H = G$, or there exists a maximal chain (1) such that $H_{i}^{\mathfrak{F}}\leq H_{i-1}$
for $ i = 1, \dots, n $.

In the case when $\mathfrak{F}$ coincides with the class $\mathfrak{N}$ of all nilpotent
groups, every $\mathfrak{N}$-subnormal subgroup is
subnormal, the converse is not true in general. However, in soluble groups these concepts are equivalent.

  Another generalization of subnormal subgroups was proposed by O.~Kegel [21]. We give it according to [20, p.~236].

A subgroup $H$ is called {\it$ {\rm K}$-$\mathfrak{F}$-subnormal} in $G$ $($which is denoted by
$H$ ${\rm K}$-$\mathfrak{F}$-${\rm sn}$ $G)$ if
there is a chain of subgroups (1) such
either $ H_ {i-1}\trianglelefteq H_{i}$, or $H_{i}^{\mathfrak{F}}\leq H_{i-1}$ for $i = 1, \dots, n$.

 Note that a subnormal subgroup is $\rm{K}$-$\mathfrak{F}$-subnormal in any group, the converse is not always true.  For the case  $\mathfrak{F}=\mathfrak{N}$ the concepts of subnormal and
$\rm{K}$-$\mathfrak{N}$-subnormal subgroups are equivalent.
If $\mathfrak{F}$ coincides with the class $\mathfrak{U}$ of all supersoluble groups, then the concept of $\mathbb{P}$-subnormal subgroup is equivalent to the concept of $\mathfrak{U}$-subnormal and $\rm {K}$-$\mathfrak{U}$-subnormal subgroup in the class of all  soluble groups.
In an arbitrary group, every $\mathfrak{U}$-subnormal ($\rm{K}$-$\mathfrak{U}$-subnormal)
subgroup is $\mathbb{P}$-subnormal ($\rm{K}$-$\mathbb{P}$-subnormal subgroup, respectively), but the converse fails in general.

The monograph [21] reflects the results of many papers in which the properties of $\mathfrak{F}$-subnormal, ${\rm K}$-$\mathfrak{F}$-subnormal subgroups and their applications were studied.

In [22] 
consideration of the following general problem was started.
{\it Let $\mathfrak{F}$ be a non-empty formation.
How $\mathfrak{F}$-subnormal $(\rm{K}$-$\mathfrak{F}$-subnormal$)$ Sylow subgroups
influence on the structure of the whole group.}
The classes ${\rm W}_{\pi}\mathfrak{F}$ and $\overline{{\rm W}}_{\pi}\mathfrak{F}$ were investigated in [23]; where ${\rm W}_{\pi}\mathfrak{F}$ ($\overline{{\rm W}}_{\pi}\mathfrak{F}$) is the class of all groups $G$, for which 1 and all Sylow $p$-subgroups are $\mathfrak{F}$-subnormal (respectively $\rm {K}$-$\mathfrak{F}$-subnormal) in $G$ for every  $p\in\pi\cap\pi(G)$.
The classes ${\rm W}\mathfrak{F}$ and $\overline {{\rm W}}\mathfrak{F}$ ($\pi$ coincides with the set of all primes) were studied in [24-27]. An interesting generalization of classes ${\rm W}_{\pi}\mathfrak{F}$ and $\overline{{\rm W}}_{\pi}\mathfrak{F}$ was considered in
[28].

\medskip

{\bf Definition 1} [29]. {\it Let $\mathfrak{F}$ be a non-empty formation.  A subgroup $H$ of a group $G$
is called strongly ${\rm K}$-$\mathfrak{F}$-subnormal in $G$,
 if  $N_{G}(H)$ is a $\mathfrak{F}$-subnormal subgroup in $G$.
}

\medskip

Note that a subgroup is normal in its normalizer. Therefore every strongly ${\rm K}$-$\mathfrak{F}$-subnormal subgroup is ${\rm K}$-$\mathfrak{F}$-subnormal in any group. The converse is not true.
Let $S$ be a symmetric group of degree 3.
By [29, theorem B. 10.9] $S$ has an irreducible and faithful $S$-module  $U$ over the field $\mathrm{F}_{7}$ of 7 elements.
Consider the semidirect product $G = [U]S$. The group $G$ is not supersoluble, because $S$ is non-abelian.
 Since $G/U$ is supersoluble, we see that $H = UQ$ is ${\rm K}$-$\mathfrak{U}$-subnormal subgroup of $G$,
where $Q$ is a Sylow $3$-subgroup of $G$ that is contained in $S$. Since $H$ is supersoluble, we deduced that $Q$ is $ {\rm K}$-$\mathfrak{U}$-subnormal in $G$. Note that the subgroup $Q$ is not strongly ${\rm K}$-$\mathfrak{U}$-subnormal in $G$. This follows from the fact that $N_{G}(Q)=S$, but $S$ is not normal and not $\mathfrak{U}$-subnormal in $G$.

\medskip

{\bf Definition 2} [29]. {\it Given a set of primes $\pi$ and a non-empty formation $\mathfrak{F}$. Introduce the following class of groups:
 ${\rm w}^{*}_{\pi}\mathfrak{F}$ is the class of all groups $G$, for which $\pi(G)\subseteq \pi(\mathfrak{F})$ and all its Sylow $q$-subgroups  are strongly $\mathfrak{F}$-subnormal in $G$ for every $q\in\pi\cap \pi(G)$.
}

When $\pi = \mathbb{P}$ is the set of all primes, we denote
${\rm w}^{*}_{\mathbb{P}}\mathfrak{F}={\rm w}^{*}\mathfrak{F}$. If  $\pi(G)\subseteq \pi(\mathfrak{F})$ and $\pi\cap \pi(G)=\emptyset$,
 then  $N_{G}(1)=G$ is $\mathfrak{F}$-subnormal in $G$ and $G\in {\rm w}^{*}_{\pi}\mathfrak{F}$.

\medskip

{\bf Problem.} {\it Let $\mathfrak{F}$ be a hereditary saturated formation
and $\pi$ be some set of primes.

$(1)$ Investigate how the properties of the class ${\rm w}^{*}_{\pi}\mathfrak{F}$ depend on the corresponding properties of $\mathfrak{F}$.
In particular, find conditions under which the class ${\rm w}^{*}_{\pi}\mathfrak{F}$ is also a saturated formation;

$(2)$ Describe  $\mathfrak{F}$ for which ${\rm w}^{*}_{\pi}\mathfrak{F} = \mathfrak{F}$.
}

\medskip

This paper is devoted studying for some cases of this problem.

\bigskip

{\parindent=0mm
\textbf{\bf {\large 1. Preliminary results}}
}

\bigskip

We use standard notation and definitions. The appropriate information on groups theory and formations theory can be found in monographs [19], [20] and [30].
We recall some concepts significant in the paper.

By  $\mathbb{P}$ we denote the set of all primes. If $\pi\subseteq\mathbb{P}$, then $\pi'=\mathbb{P}\setminus\pi$.
Let  $G$ be a group and $p$ be a prime.
We denote by $|G|$ the order of $G$;
by $\pi(G)$, the set of all prime divisors of $|G|$;
by $O_{p}(G)$, the largest normal $p$-subgroup of $G$;
by $O_{\pi}(G)$, the largest normal $\pi$-subgroup of $G$;
by ${\rm Syl}_{p}(G)$, the set of all Sylow $p$-subgroups of $G$;
by ${\rm Syl}(G)$, the set of all Sylow subgroups of $G$;
by $F(G)$, the Fitting subgroup of $G$, which is the largest normal nilpotent  subgroup of $G$;
by $F_p(G)$, the $p$-nilpotent radical of $G$, which  is the largest normal $p$-nilpotent  subgroup of $G$;
by $Z_p$, the cyclic group of order $p$;
by $1$, the  identity subgroup (group).

By  $l_{p}(G)$ we denote the $p$-length of the $p$-soluble group $G$;
an arithmetic length of the soluble group $G$ is $al(G) = \mathrm{Max}\,l_{p}(G)$, where $p$ runs through all primes $p\in\pi(G)$;
$\mathfrak{L}_{a}(n)$ is the class of all soluble groups $G$ with  $al(G)\leq n$;
$\mathfrak{L}_{a}(1)$ is the class of all soluble groups $G$ with $al(G)\leq 1$.

In the next lemma, the some familiar properties of Sylow subgroups are collected.

\medskip

{\bf Lemma~1.1.} 
{\it Let $G$ be a group and $p \in \mathbb{P}$. Then the following statements are true.

$(1)$ If $P\in {\rm Syl}_{p}(G)$ and $N \trianglelefteq G$, then
$P\cap N\in{\rm Syl}_{p}(N)$ and
$PN/N \in{\rm Syl}_{p}(G/N)$, moreover $N_{G/N}(PN/N)=N_{G}(P)N/N$.

$(2)$ If $H/N\in {\rm Syl}_{p}(G/N)$ and $N \trianglelefteq G$, then $H/N=PN/N$ for some $P\in {\rm Syl}_{p}(G)$.

$(3)$ If $P \in {\rm Syl}(G)$ and $N_{i} \trianglelefteq G$, $i=1, 2$, then

\centerline{
$P \cap N_{1}N_{2} = (P \cap N_{1})(P \cap N_{2})$ and $PN_{1} \cap PN_{2} = P(N_{1} \cap N_{2})$.}

$(4)$ If $\pi(G)=\{p_1, \ldots, p_r\}$ and $P_i \in {\rm Syl}_{p_i}(G)$ for
$i=1, \ldots, r$, then $G= \langle P_1, \ldots, P_r \rangle$. 
}

\medskip

{\bf Lemma~1.2} [30, lemma~A.1.2] 
{\it Let $U$, $V$ and $W$ be subgroups of $G$. Then the following statements are equivalent:

$(1)$ $U\cap VW=(U\cap V)(U\cap W);$

$(2)$ $UV\cap UW=U(V\cap W)$.
}

\medskip

{\bf Proposition~1.3.}
{\it Let $G$ be a group, $P \in {\rm Syl}(G)$ and $N_{i} \trianglelefteq G$, $i=1, 2$. Then

\centerline{
$N_{G}(P)\cap N_{1}N_{2}=(N_{G}(P)\cap N_{1})(N_{G}(P)\cap N_{2})$ and}
\centerline{
$N_{G}(P)N_{1}\cap N_{G}(P)N_{2}=N_{G}(P)(N_{1}\cap N_{2})$.
}
}
\medskip

{\sc Proof.} We proceed  by induction on $|G|$. Let $N_{1}$ and $N_{2}$ be normal subgroups of $G$ and $P \in {\rm Syl}(G)$. If $N_{1}\cap N_{2}\neq 1$, then there exist a minimal normal subgroup $N$ of $G$, contained in $N_{1}\cap N_{2}$. By induction
	
\centerline{
		$N_{G/N}(PN/N)\cap N_{1}/N\cdot N_{2}/N=(N_{G/N}(PN/N)\cap N_{1}/N)(N_{G/N}(PN/N)\cap N_{2}/N)$.
}

	By Lemma 1.1(1) $N_{G/N}(PN/N)=N_{G}(P)N/N$. By the Dedekind identity, we have
		$N_{G}(P)N/N\cap N_{1}N_{2}/N=(N_{G}(P)\cap N_{1}N_{2})N/N$ and
		$N_{G}(P)N/N\cap N_{i}/N=(N_{G}(P)\cap N_{i})N/N$ for $i=1,2$.
	
Then $N_{G}(P)\cap N_{1}N_{2}=N_{G}(P)\cap (N_{G}(P)N\cap N_{1}N_{2})=N_{G}(P)\cap (N_{G}(P)\cap N_{1})N\cdot (N_{G}(P)\cap N_{2})N=(N_{G}(P)
\cap N_{1})(N_{G}(P)\cap N_{2})(N_{G}(P)\cap N)=(N_{G}(P)\cap N_{1})(N_{G}(P)\cap N_{2})$.
	
	Let $N_{1}\cap N_{2}=1$. Let $T=N_{G}(P)N_{1}\cap N_{G}(P)N_{2}$. Since $PN_{i}\unlhd N_{G}(P)N_{i}$, $i=1,2$, we have $PN_{1}\cap PN_{2}\unlhd T$.
From $N_{1}\cap N_{2}=1$ and lemma 1.1(3) 
follows that $PN_{1}\cap PN_{2}=P(N_{1}\cap N_{2})=P$. Therefore $P\unlhd T$ and $T=N_{G}(P)$. Then
$N_{G}(P)(N_{1}\cap N_{2})=N_{G}(P)=N_{G}(P)N_{1}\cap N_{G}(P)N_{2}$.
	By lemma 1.2	
			$N_{G}(P)\cap N_{1}N_{2}=(N_{G}(P)\cap N_{1})(N_{G}(P)\cap N_{2})$. $\Box$
	
\medskip

{\bf Lemma~1.4.} [19,~lemma~3.9]. 
{\it If $H / K$ is a chief factor of a group $G$  and $p\in \pi(H/K)$,
then $G/C_{G}(H/K)$ does not contain nonidentity normal $p$-subgroups, and
$F_{p}(G)\leq C_{G}(H/K)$.
}

\medskip

Let  $\mathfrak{F}$ be a class of groups.
By $\pi(\mathfrak{F})$
we denote the set of all prime divisors of orders of groups belonging to $\mathfrak{F}$; $\mathfrak{F}_\pi$ is the class of all
$\pi$-groups belonging to $\mathfrak{F}$;
$\mathfrak{F}_{p} = \mathfrak{F}_{\pi}$ for $\pi = \{p\}$.

We will use the following notation:
$\mathfrak{G}$ is the class of all groups,
$\mathfrak{S}$ is the class of all soluble groups,
$\mathfrak{N}$ is the class of all nilpotent groups,
$\mathfrak{N}^{2}$ is the class of all metanilpotent groups,
$\mathfrak{NA}$ is the class of all groups $G$ with the nilpotent commutator subgroup $G'$.

{\it A minimal non-$\mathfrak{F}$-group} is
a group $G$ such that $G \not\in \mathfrak{F}$,
and any proper subgroup of $G$ belongs to $\mathfrak{F}$.
A minimal non-$\mathfrak{N}$-group
 is called a {\it Schmidt group}.

A class of groups $\mathfrak{F}$ is called a {\it formation}, if
1) $\mathfrak{F}$ is a homomorph, i.e., from $G \in \mathfrak{F}$
and $N \trianglelefteq G$ it follows that $G/N \in \mathfrak{F}$ and
2) from $N_i \trianglelefteq G$ and $G/N_i \in \mathfrak{F}$ $(i=1, 2)$
it ensues that $G/N_1 \cap N_2 \in \mathfrak{F}$.

A formation $\mathfrak{F}$ is called {\it saturated}, if
from $G/\Phi(G)\in \mathfrak{F}$ it follows that
$G \in \mathfrak{F}$. A formation $\mathfrak{F}$ is called {\it hereditary}
if, together with each group, $\mathfrak{F}$ contains all its subgroups.
By symbol $G^{\mathfrak{F}}$ denotes the {\it $\mathfrak{F}$-residual} of $G$; i.e., the least normal subgroup of $G$ for which $G/G^{\mathfrak{F}} \in \mathfrak{F}$.

A function $f : \mathbb{P} \to \{$formations$\}$ is called a {\it local screen}.
By $LF(f)$ we denote the class of all groups $G$ with $G/C_{G}(H/K) \in
f(p)$ for each chief factor $H/K$ and each $p \in \pi(H/K)$.
A formation $\mathfrak{F}$ is called  {\it local}, if there exists a local screen $f$
with $\mathfrak{F} = LF(f)$.

A screen $f$ of a formation $\mathfrak{F}$  is called {\it inner} if
$f(p) \subseteq \mathfrak{F}$ for each prime $p$.
 An inner screen $f$ of $\mathfrak{F}$ is called the {\it maximal inner}
if, for its every inner screen $h$, we have $h(p) \subseteq f(p)$ for every prime $p$.

\medskip

{\bf Lemma~1.5} [19, lemma~4.5]. 
{\it Let $\mathfrak{F} = LF(f)$.
A group  $G$ belongs to $\mathfrak{F}$ if and only if $G/F_{p}(G)\in f(p)$
for each $p \in \pi(G)$.
}

\medskip

 We give some knows properties of $\mathfrak{F}$-subnormal and $\rm{K}$-$\mathfrak{F}$-subnormal subgroups.

\medskip

{\bf Lemma~1.6.}
{\it Let  $\mathfrak{F}$
be a non-empty formation, $H$ and $K$ are subgroups of a group $G$,
and $N \trianglelefteq G$.

$(1)$ If $H$ $\mathfrak{F}$-${\rm sn}$ $G$ $(H$ ${\rm K}$-$\mathfrak{F}$-${\rm sn}$ $G)$
then $HN/N$ $\mathfrak{F}$-${\rm sn}$ $G/N$ $(HN/N$ ${\rm
K}$-$\mathfrak{F}$-${\rm sn}$ $G/N).$

$(2)$ If  $N\leq H$ and $H/N$ $\mathfrak{F}$-${\rm sn}$ $G/N$
$(H/N$ ${\rm K}$-$\mathfrak{F}$-${\rm sn}$ $G/N)$ then  $H$ $\mathfrak{F}$-${\rm sn}$ $G$
$(H$ ${\rm K}$-$\mathfrak{F}$-${\rm sn}$ $G).$

$(3)$ If  $H$ $\mathfrak{F}$-${\rm sn}$ $G$ $(H$ ${\rm K}$-$\mathfrak{F}$-${\rm sn}$ $G)$
then  $HN$ $\mathfrak{F}$-${\rm sn}$ $G$ $(HN$ ${\rm
K}$-$\mathfrak{F}$-${\rm sn}$ $G).$

$(4)$ If  $H$ $\mathfrak{F}$-${\rm sn}$ $K$ $(H$ ${\rm K}$-$\mathfrak{F}$-${\rm sn}$ $K)$
and  $K$ $\mathfrak{F}$-${\rm sn}$ $G$ $(K$ ${\rm K}$-$\mathfrak{F}$-${\rm sn}$ $G)$
then  $H$ $\mathfrak{F}$-${\rm sn}$ $G$ $(H$ ${\rm K}$-$\mathfrak{F}$-${\rm
sn}$ $G).$

$(5)$ If all composition factors of $G$ belong to $\mathfrak{F}$ then every subnormal subgroup of $G$  is $\mathfrak{F}$-subnormal$.$

$(6)$ Let $p$ be a prime and let $G$ be a $p$-group. If $Z_{p}\in \mathfrak{F}$
then all subgroups of $G$ are $\mathfrak{F}$-subnormal.
}

\medskip

{\bf Lemma~1.7.}
{\it Let $\mathfrak{F}$
be a non-empty hereditary formation, $H \leq G$ and $M \leq G$.

$(1)$ If  $H$ $\mathfrak{F}$-${\rm sn}$ $G$ $(H$ ${\rm K}$-$\mathfrak{F}$-${\rm sn}$ $G)$
then  $H \cap M$ $\mathfrak{F}$-${\rm sn}$ $M$ $(H \cap M $ ${\rm
K}$-$\mathfrak{F}$-${\rm sn}$ $M).$

$(2)$ If  $H$ $\mathfrak{F}$-${\rm sn}$ $G$ and $M$ $\mathfrak{F}$-${\rm sn}$ $G$
$(H$ ${\rm K}$-$\mathfrak{F}$-${\rm sn}$ $G$ and $M$ ${\rm K}$-$\mathfrak{F}$-${\rm sn}$ $G)$
then  $H \cap M $ $\mathfrak{F}$-${\rm sn}$ $G$ $(H \cap M$ ${\rm
K}$-$\mathfrak{F}$-${\rm sn}$ $G).$

$(3)$ If $G^{\mathfrak{F}} \leq H$ then $H$ $\mathfrak{F}$-${\rm sn}$ $G$
$(H$ ${\rm K}$-$\mathfrak{F}$-${\rm sn}$ $G).$

$(4)$ If $H$ $\mathfrak{F}$-${\rm sn}$ $G$ $(H$ ${\rm K}$-$\mathfrak{F}$-${\rm sn}$ $G)$
then $H^{x}$ $\mathfrak{F}$-${\rm sn}$ $G$ $(H^{x}$ ${\rm K}$-$\mathfrak{F}$-${\rm sn}$ $G)$
for any $x \in G$.
}

\bigskip

{\parindent=0mm
\textbf{\bf {\large 2. Properties of the Class ${\rm w}^{*}_{\pi}\mathfrak{F}$}}
}

\bigskip

Recall that the class of groups ${\rm w}^{*}_{\pi}\mathfrak{F}$ is defined as follows:

${\rm w}^{*}_{\pi}\mathfrak{F}=(G$ $|$  $\pi(G)\subseteq \pi(\mathfrak{F})$ and every Sylow $q$-subgroup of $G$ is strongly $\mathfrak{F}$-subnormal in $G$, where $q\in\pi\cap \pi(G))$.

The following example shows that ${\rm w}^{*}_{\pi}\mathfrak{F}\not=\mathfrak{F}$ in the general case.

\medskip

 {\bf Example 2.1.} Let $\mathfrak{F}=\mathfrak{N}^{3}$ be the formation of all soluble groups whose nilpotent length is $\leq 3$.
 Take  the symmetric group $S_{4}=M$ of degree 4.
 By [30, theorem B. 10.9] there exists an irreducible and faithful $M$-module  $U$ over the field $\mathrm{F}_{3}$ of 3 elements.
Consider the semidirect product $G = [U]M$.
 Note that the nilpotent length of $G$ is 4 and $\pi(G)= \{2,3\}$. Since $S$ is a minimal non-$\mathfrak{N}^{2}$-subgroup, we deduced that
 $G$ is minimal non-$\mathfrak{N}^{3}$-group.
 It is easy to see that the normalizers of its Sylow subgroups are $\mathfrak{F}$-subnormal subgroups in $G$, but $G$
 does not belong to $\mathfrak{F}$.

\medskip

{\bf Definition 2.2.} {\it A class of groups $\mathfrak{F}$ is called $S_{H}$-closed, if from $G\in\mathfrak{F}$ it follows that every Hall subgroup of $G$ belongs to $\mathfrak{F}$.
}

\medskip

{\bf Proposition 2.3.}
{\it Let  $\mathfrak{F}$ be a non-empty formation and $\pi \subseteq \mathbb{P}$.

$(1)$ If  $\pi_1$ is a set of primes and $\pi\subseteq\pi_1$ then
${\rm w}^{*}_{{\pi}_1}\mathfrak{F} \subseteq {\rm w}^{*}_{\pi}\mathfrak{F}.$

$(2)$ $\mathfrak{N}_{\pi \cap \pi(\mathfrak{F})} \subseteq {\rm w}^{*}_{\pi}\mathfrak{F}.$

$(3)$ ${\rm w}^{*}_{\pi}\mathfrak{F} = {\rm w}^{*}_{\pi \cap
\pi(\mathfrak{F})}\mathfrak{F}.$

$(4)$ ${\rm w}^{*}_{\pi}\mathfrak{F}$ is a homomorph.

$(5)$ If a formation $\mathfrak{F}_{1}\subseteq\mathfrak{F}$ then ${\rm w}^{*}_{\pi}\mathfrak{F}_{1}\subseteq{\rm w}^{*}_{\pi}\mathfrak{F}$.

}

\medskip

{\sc Proof.}
(1): Let  $G\in {\rm w}^{*}_{{\pi}_1}\mathfrak{F}$, $q\in \pi\cap \pi(G)$ and $Q$ be any Sylow $q$-subgroup of $G$. Since  $q\in \pi_1\cap \pi(G)$, we have
$N_{G}(Q)\ \mathfrak{F}$-${\rm sn}$ $G$. Hence
${\rm w}^{*}_{{\pi}_1}\mathfrak{F} \subseteq {\rm w}^{*}_{\pi}\mathfrak{F}$.

(2): Let $G \in \mathfrak{N}_{\pi \cap \pi(\mathfrak{F})}$. Then $\pi(G)\subseteq (\pi \cap \pi(\mathfrak{F}))\subseteq\pi(\mathfrak{F})$. Since $N_{G}(P)=G$ for every $P\in{\rm Syl}(G)$, by definition~1 it follows that $G \in{\rm w}^{*}_{\pi}\mathfrak{F}$.

(3): From (1) it follows that
${\rm w}^{*}_{\pi}\mathfrak{F} \subseteq {{\rm w}^{*}}_{\pi \cap \pi(\mathfrak{F})}\mathfrak{F}$.

Let $G \in {{\rm w}^{*}}_{\pi \cap \pi(\mathfrak{F})}\mathfrak{F}$. Since
$\pi(G) \subseteq \pi(\mathfrak{F})$, we have
$\pi \cap \pi(\mathfrak{F}) \cap \pi(G) = \pi \cap \pi(G)$. Consequently,
if $q \in \pi \cap \pi(G)$, then in $G$ the normalizer of every  Sylow $q$-subgroup is $\mathfrak{F}$-subnormal. So $G \in {\rm w}^{*}_{\pi}\mathfrak{F}$ and
${\rm w}^{*}_{\pi}\mathfrak{F} = {{\rm w}^{*}}_{\pi \cap \pi(\mathfrak{F})}\mathfrak{F}$.

(4): To prove that ${\rm w}^{*}_{\pi}\mathfrak{F}$ is a homomorph,
let  $G \in {\rm w}^{*}_{\pi}\mathfrak{F}$, $N \trianglelefteq G$ and $p \in \pi\cap\pi(G/N)$.
Consider $H/N \in {\rm Syl}_{p}(G/N)$.
By Lemma 1.1(2) $H/N = PN/N$ for some Sylow $p$-subgroup $P$ of $G$.
From  $G \in {\rm w}^{*}_{\pi}\mathfrak{F}$  it follows that $N_{G}(P)$ $\mathfrak{F}$-${\rm sn}$ $G$.
Then by Lemma~1.1(1) and Lemma 1.6(1)
$N_{G/N}(H/N) = N_{G}(P)N/N$ $\mathfrak{F}$-${\rm sn}$ $G/N$. From here and $\pi(G/N)\subseteq\pi(G)\subseteq\pi(\mathfrak{F})$ we have that $G/N \in {\rm w}^{*}_{\pi}\mathfrak{F}$. So
${\rm w}^{*}_{\pi}\mathfrak{F}$ is a homomorph.

(5): Let  $G\in{\rm w}^{*}_{\pi}\mathfrak{F}_{1}$. Then $\pi(G)\subseteq\pi(\mathfrak{F}_{1})\subseteq\pi(\mathfrak{F})$.
From $q\in\pi\cap \pi(G)$ it follows that every $Q\in{\rm Syl}_{q}(G)$ is strongly ${\rm K}$-$\mathfrak{F}_{1}$-subnormal in $G$. If $N_{G}(Q)= G$, then $N_{G}(Q)$ $\mathfrak{F}$-${\rm sn}$ $G$. Suppose that a maximal chain of subgroups
$N_{G}(Q) = H_{0} < H_{1} < \dots < H_{n} = G$ exists and $H_{i}^{\mathfrak{F}_{1}} \leq H_{i-1}$
for $i = 1, \dots, n$.
 From $H_{i}/H_{i}^{\mathfrak{F}_{1}}\in\mathfrak{F}_{1}\subseteq\mathfrak{F}$  we have $H_{i}^{\mathfrak{F}}\subseteq H_{i}^{\mathfrak{F}_{1}}\leq H_{i-1}$. Hence $N_{G}(Q)$ $\mathfrak{F}$-${\rm sn}$ $G$. So ${\rm w}^{*}_{\pi}\mathfrak{F}_{1}\subseteq{\rm w}^{*}_{\pi}\mathfrak{F}$.
$\Box$

\medskip

{\bf Theorem 2.4.}
{\it Let  $\mathfrak{F}$ be a non-empty hereditary formation and $\pi \subseteq \mathbb{P}$.
Then 

$(1)$ $\mathfrak{F} \subseteq {\rm w}^{*}\mathfrak{F}\subseteq {\rm w}^{*}_{\pi}\mathfrak{F}$,

$(2)$ ${\rm w}^{*}_{\pi}\mathfrak{F}$ is an $S_{H}$-closed formation,

$(3)$ ${\rm w}^{*}_{\pi}\mathfrak{F}={\rm w}^{*}_{\pi}({\rm w}^{*}_{\pi}\mathfrak{F})$.

}

\medskip

{\sc Proof.}
(1): From Lemma~1.7(3) it follows that $\mathfrak{F} \subseteq {\rm w}^{*}\mathfrak{F}$.
From $\pi \subseteq \mathbb{P}$ and Proposition~2.3(1) we conclude that
${\rm w}^{*}\mathfrak{F}\subseteq {\rm w}^{*}_{\pi}\mathfrak{F}$.

(2): To prove $S_{H}$-closure of ${\rm w}^{*}_{\pi}\mathfrak{F}$,
let  $G \in {\rm w}^{*}_{\pi}\mathfrak{F}$
and let $H$ be a Hall subgroup of $G$. Then $\pi(H)\subseteq\pi(G)\subseteq\pi(\mathfrak{F})$. Let $q \in \pi\cap\pi(H)$ and $S$ be a Sylow $q$-subgroup
of $H$. Since $S\in {\rm Syl}_{q}(G)$, we have $N_{G}(S)$ $\mathfrak{F}$-${\rm sn}$ $G$. By Lemma~1.7(1) $N_{H}(S)=(N_{G}(S)\cap H)$ $\mathfrak{F}$-${\rm sn}$ $H$. Therefore $H \in {\rm w}^{*}_{\pi}\mathfrak{F}$.

By Proposition~2.3(4) ${\rm w}^{*}_{\pi}\mathfrak{F}$ is a homomorph.

Let us proved that ${\rm w}^{*}_{\pi}\mathfrak{F}$ is closed under subdirect products. Suppose that is false, and let $G$ be a counterexample with $|G|$ as small as possible. Then there exists a subgroup $N_{i}\unlhd G$  such that  $G/N_{i} \in {{\rm w}^{*}_{\pi}}\mathfrak{F}$, $i = 1, 2$, but
 $G/N_{1} \cap N_{2} \notin {{\rm w}^{*}_{\pi}}\mathfrak{F}$.
We note that from $\pi(G/N_{i})\subseteq\pi(\mathfrak{F})$, $i=1, 2$, it follows that $\pi(G/N_{1}\cap N_{2})\subseteq\pi(\mathfrak{F})$.
By the choice of $G$ we can assume that $N_{1} \cap N_{2} = 1$.
Let $p \in \pi\cap\pi(G)$ and $R\in{\rm Syl}_{p}(G)$.
Since  $RN_{i}/N_{i}$ is a Sylow $p$-subgroup
of $G/N_{i}$ and $G/N_{i} \in {{\rm w}^{*}_{\pi}}\mathfrak{F}$,
we have $N_{G/N_{i}}(RN_{i}/N_{i})$ $\mathfrak{F}$-${\rm sn}$ $G/N_{i}$, $i = 1, 2$.
By Lemmas 1.1(1) and 1.6(2) $N_{G}(R)N_{i}$ $\mathfrak{F}$-${\rm sn}$ $G$,
$i = 1, 2$. From Lemma 1.7(2) it follows $N_{G}(R)N_{1}\cap N_{G}(R)N_{2}$ $\mathfrak{F}$-${\rm sn}$ $G$.
From Proposition 1.3 we conclude that
$N_{G}(R)N_{1}\cap N_{G}(R)N_{2}= N_{G}(R)(N_{1} \cap N_{2}) = N_{G}(R)$ $\mathfrak{F}$-${\rm sn}$ $G$.
We have the contradiction to the choice of $G$. So ${\rm w}^{*}_{\pi}\mathfrak{F}$
is closed under subdirect products.

(3): Denote $\mathfrak{X}={{\rm w}^{*}_{\pi}}\mathfrak{F}$.  Let $G \in \mathfrak{X}$. Then $\pi(G)\subseteq\pi(\mathfrak{F})$. By (1) we have that $\mathfrak{F}\subseteq\mathfrak{X}$. Therefore $\pi(G)\subseteq\pi(\mathfrak{X})$.
Let  $q \in \pi\cap\pi(G)$ and $Q\in{\rm Syl}_{q}(G)$. From $G \in \mathfrak{X}$ it  follows that $N_{G}(Q)$ $\mathfrak{F}$-${\rm sn}$ $G$. Assume that $N_{G}(Q)\not=G$. Then there is a maximal chain of subgroups $N_{G}(Q) = H_{0} < H_{1} < \dots < H_{n} = G$ such that $H_{i}^{\mathfrak{F}} \leq H_{i-1}$
for $i = 1, \dots, n$. By (2) $\mathfrak{X}$  is a formation. Therefore from $H_{i}/H_{i}^{\mathfrak{F}}\in\mathfrak{F}\subseteq\mathfrak{X}$ it follows that $H_{i}^{\mathfrak{X}}\leq H_{i}^{\mathfrak{F}}\leq H_{i-1}$. This means that $N_{G}(Q)$ $\mathfrak{X}$-${\rm sn}$ $G$. If $N_{G}(Q)=G$, then $N_{G}(Q)$ $\mathfrak{X}$-${\rm sn}$ $G$. So $G\in{{\rm w}^{*}_{\pi}}\mathfrak{X}$ and $\mathfrak{X}\subseteq{{\rm w}^{*}_{\pi}}\mathfrak{X}$ is proved.

Suppose that $\mathfrak{X}\not={{\rm w}^{*}_{\pi}}\mathfrak{X}$.
Let  $G$ be the group of minimal order in
${{\rm w}^{*}_{\pi}}\mathfrak{X}\textbackslash \mathfrak{X}$. Then $\pi(G)\subseteq\pi(\mathfrak{X})\subseteq\pi(\mathfrak{F})$. Since $G\not\in\mathfrak{X}$, there exists $P\in{\rm Syl}_{p}(G)$ such that $p\in\pi\cap\pi(G)$ and $N_{G}(P)$  is not $\mathfrak{F}$-subnormal in $G$.
We note that $N_{G}(P)$ $\mathfrak{X}$-${\rm sn}$ $G$.
Then $N_{G}(P)\not=G$ and there exists a maximal chain of subgroups $N_{G}(P) = H_{0} < H_{1} < \dots < H_{n-1} < H_{n}= G$ such that $H_{i}^{\mathfrak{X}} \leq H_{i-1}$
for $i = 1, \dots, n$. Since $N_{G}(P) = N_{H_{i}}(P)$, $N_{H_{i}}(P)H_{i}^{\mathfrak{X}} \leq H_{i-1}$ and $H_{i}/H_{i}^{\mathfrak{X}}\in\mathfrak{X}$, we have $N_{H_{i}}(P)H_{i}^{\mathfrak{X}}/H_{i}^{\mathfrak{X}}=
N_{H_{i}/H_{i}^{\mathfrak{X}}}(PH_{i}^{\mathfrak{X}}/H_{i}^{\mathfrak{X}})$ $\mathfrak{F}$-${\rm sn}$ $H_{i}/H_{i}^{\mathfrak{X}}$. By Lemma~1.6(2) $N_{H_{i}}(P)H_{i}^{\mathfrak{X}}$ $\mathfrak{F}$-${\rm sn}$ $H_{i}$ for $i = 1, \dots, n$. Therefore $H_{n}^{\mathfrak{X}}=G^{\mathfrak{X}}\not\subseteq N_{G}(P)$. From the maximality of $N_{G}(P)$ in $H_{1}$ it follows that $N_{G}(P)$ $\mathfrak{F}$-${\rm sn}$ $H_{1}$. So $n\not=1$. Suppose that $n=2$. Then by Lemma~1.7(1) $N_{G}(P)=N_{G}(P)\cap N_{G}(P)H_{2}^{\mathfrak{X}}$ $\mathfrak{F}$-${\rm sn}$ $N_{G}(P)H_{2}^{\mathfrak{X}}$. From $N_{G}(P)H_{2}^{\mathfrak{X}}$ $\mathfrak{F}$-${\rm sn}$ $H_{2}$ we conclude that $N_{G}(P)$ $\mathfrak{F}$-${\rm sn}$ $H_{2}=G$. This is the contradiction with the choice of $G$. So, we can assume that $n\geq 3$ and $N_{G}(P)$ $\mathfrak{F}$-${\rm sn}$ $H_{n-1}$. Since  $N_{G}(P)H_{n}^{\mathfrak{X}}\leq H_{n-1}$, by Lemma~1.7(1) we have $N_{G}(P)=N_{G}(P)\cap N_{G}(P)H_{n}^{\mathfrak{X}}$ $\mathfrak{F}$-${\rm sn}$ $N_{G}(P)H_{n}^{\mathfrak{X}}$. From $N_{G}(P)H_{n}^{\mathfrak{X}}$ $\mathfrak{F}$-${\rm sn}$ $G$ it follows that $N_{G}(P)$ $\mathfrak{F}$-${\rm sn}$ $G$. This contradicts the choice of $G$. So $\mathfrak{X}={{\rm w}^{*}_{\pi}}\mathfrak{X}$.
$\Box$

\bigskip

{\parindent=0mm
\textbf{\bf {\large 3. Formations $\mathfrak{F}$ for which ${\rm
w}^{*}_{\pi}\mathfrak{F}=\mathfrak{F}$}}
}


\bigskip

This section focuses on (2) of Problem.

\medskip

{\bf Lemma~3.1.}
{\it $(1)$ The class $\mathfrak{L}_{a}(1)$ is a hereditary saturated Fitting formation.

$(2)$ Let  $G$ be a soluble group, $\Phi(G)=1$. $G$ is a minimal non-$\mathfrak{L}_{a}(1)$-group if and only if the following statements hold:

$1)$ $|G|=p^{\alpha}q^{\beta}$, $l_{p}(G)=1$, $l_{q}(G)=2$, $l(G)=3$;

$2)$ $G$ has precisely three conjugate classes of maximal subgroups, whose representatives have the following structure:
$G_{p}\leftthreetimes G_{q}^{*}$, the Schmidt group, $F(G)\leftthreetimes G_{p}$ and $G_{q}\leftthreetimes \Phi(G_{p})$,
where $G_{q}=F(G)\leftthreetimes G_{q}^{*}$.} 

\medskip

{\sc Proof.}  (1): The statement follows directly from the fact that $\mathfrak{L}_{a}(1)=\cap\mathfrak{G}_{p'}\mathfrak{N}_{p}\mathfrak{G}_{p'}$ for all $p\in\mathbb{P}$.

(2): The statement is Lemma 4.1 in [31]. $\Box$

\medskip

{\bf Lemma~3.2.}
{\it Let  $G$ be a biprimary group and let $G\in\mathfrak{L}_{a}(1)$. Then  $G$ is metanilpotent.}

{\sc Proof.} Let  $G$ be a counterexample of minimal order to the statement of the lemma. Since $\mathfrak{N}^{2}$ is a hereditary saturated formation, the group $G=NM$, where  $N$ is a unique minimal normal subgroup of $G$ and $M$ is a maximal subgroup of $G$, moreover, $N$ is an abelian $p$-group,  $p$ is some prime, $M$ is a Schmidt group with a normal $p$-subgroup. From $O_{p}(M)=1$ we conclude that  $p$-length of $G$ is 2. This contradicts the fact that $G\in\mathfrak{L}_{a}(1)$. $\Box$

{\bf Lemma~3.3.}
{\it Let  $\mathfrak{F}$  be a non-empty hereditary formation and let $G$ be a soluble group. If $G\in\mathfrak{L}_{a}(1)$, $G\not=N_{G}(P)$ and $N_{G}(P)\in\mathfrak{F}$ for all $P\in{\rm Syl}(G)$, then $G\in\mathfrak{F}$.}

{\sc Proof.} Let  $G$ be a counterexample of minimal order to the statement of the lemma. Let $N$ is a minimal normal subgroup of $G$. We will prove that $G/N\in\mathfrak{F}$. If $G/N\not=N_{G/N}(H/N)$ for all $H/N\in{\rm Syl}(G/N)$, then $G/N\in\mathfrak{F}$ by the choice of $G$.
If $G/N =N_{G/N}(H/N)$ for some $H/N\in{\rm Syl}_{q}(G/N)$, then $H/N=QN/N$ for some $Q\in{\rm Syl}_{q}(G)$ and $G=N_{G}(Q)N$.
Since $G$ is solvable and $G\not=N_{G}(Q)$, we conclude that $N_{G}(Q)$ is a maximal subgroup of $G$ and $N_{G}(Q)\cap N=1$.
From here $G/N\cong N_{G}(Q)\in\mathfrak{F}$. If  $K$  is a minimal normal subgroup of $G$ and $K\not= N$, then $G/K\in\mathfrak{F}$. Since $\mathfrak{F}$ is a formation,  we deduce that $G/N\cap K\cong G\in\mathfrak{F}$. This contradicts to the choice of $G$. Consequently $N$ is the unique minimal normal subgroup of $G$. Since $G$ is soluble, we conclude that $N$ is a $p$-group. From the uniqueness of $N$  it follows that $F(G)$ is a $p$-group.
By the choice of $G$ imply that $\pi(G)\geq 2$. If  $F(G)\not=P\in{\rm Syl}_{p}(G)$, then $p\in\pi(G/F(G))$. This contradicts with $G\in\mathfrak{L}_{a}(1)$. Therefore, $F(G)=P\in{\rm Syl}_{p}(G)$ and $G=N_{G}(P)$. This contradiction completes the proof of the lemma. $\Box$

\medskip

{\bf Theorem  3.4.}
{\it Let  $\mathfrak{F}$ be a hereditary saturated formation and  $\mathfrak{F}\subseteq\mathfrak{L}_{a}(1)$.
A group $G\in\mathfrak{F}$ if and only if
$\pi(G)\subseteq \pi(\mathfrak{F})$ and all its Sylow subgroups are strongly ${\rm K}$-$\mathfrak{F}$-subnormal in $G$.
}

\medskip

{\sc Proof.}
Necessity. Let $G\in\mathfrak{F}$. By Lemma~1.7(3) $N_{G}(S)$
$\mathfrak{F}$-${\rm sn}$ $G$ for any Sylow subgroup $S$ of $G$. 

Sufficiency. Let $G$ be a counterexample of minimal order and let $N$ be a minimal normal subgroup of $G$.

If  $G=N$ then $G$ is a simple group, because $N$ is the minimal normal subgroup of $G$. If  $G\cong Z_{p}$ then from $\pi(G)\subseteq\pi(\mathfrak{F})$ it follows that
  $G\in\mathfrak{F}$. This is the contradiction to the choice of $G$. Suppose $G$  is a simple non-abelian group and $p\in\pi(G)$. Let $G_{p}\in{\rm Syl}_{p}(G)$. Then $N_{G}(G_{p})\neq G$. From $G\notin\mathfrak{F}$ it follows that $G^{\mathfrak{F}}=G$. By hypothesis
 $N_{G}(G_{p})$ $\mathfrak{F}$-${\rm sn}$ $G$. Then there is a maximal  subgroup  $M$ of $G$ such that
 $N_{G}(G_{p})\subseteq M$ and $G^{\mathfrak{F}}\subseteq M$. This is the contradiction with $G^{\mathfrak{F}}=G$.

Let $N\neq G$. From (1)--(2) of Lemma~1.1, (1) of Lemma~1.6 and hypothesis we have $N_{G/N}(H/N)$ $\mathfrak{F}$-$\mathrm{sn}$ $G/N$ for all $H/N\in{\rm Syl}_{q}(G/N)$.
By the choice of $G$ we obtain that $G/N\in\mathfrak{F}$. If  $K$ is a minimal normal subgroup of
 $G$ and $K\neq N$, then $G/K\in \mathfrak{F}$. Since $\mathfrak{F}$ is a formation, we conclude that $G/N\cap K \cong G\in\mathfrak{F}$. This is the contradiction with the choice of $G$. Hence $G$ has the unique minimal normal subgroup $N$. If  $\Phi(G)\neq 1$, then from $G/\Phi(G)\in\mathfrak{F}$ and saturation $\mathfrak{F}$ it follows that $G\in\mathfrak{F}$. This contradicts our assumption. Therefore $\Phi(G)=1$. In this case
 $N=G^{\mathfrak{F}}$ and there is a maximal subgroup $M$ in $G$ such that $G=NM$.
Consider the following cases.

\medskip

1. $N$ is a non-abelian group. Let  $p\in\pi(N)$ and let $G_{p}\in {\rm Syl}_{p}(G)$. Then  $N_{G}(G_{p})\neq G$. Otherwise $G_{p}\unlhd G$ and
 $N\subseteq G_{p}$, since $N$ is the unique minimal normal subgroup of $G$. But then $N$  is an abelian group. This is contradiction with the proposition.

Consider  $N_{G}(G_{p})N$. Let  $N_{G}(G_{p})N = G$. From $N_{G}(G_{p})$ $\mathfrak{F}$-$\mathrm{sn}$ $G$ we deduce that there is a maximal subgroup $W$ of $G$  such that $N_{G}(G_{p})\subseteq W$ and $N=G^{\mathfrak{F}}\subseteq W$. So we have the contradiction $G=N_{G}(G_{p})N\subseteq W\neq G$.

Now let $N_{G}(G_{p})N\neq G$. Note that $G_{p}\cap N=N_{p}\in {\rm Syl}_{p}(N)$ and $N_{p}=G_{p}\cap N\unlhd N_{G}(G_{p})\cap N$. Since $N_{G}(G_{p})$ $\mathfrak{F}$-$\mathrm{sn}$ $G$, we see that ($N_{G}(G_{p})\cap N$) $\mathfrak{F}$-$\mathrm{sn}$ $N$ by Lemma~1.7(1). Since $N$ is a minimal normal subgroup of $G$, we have either $N^{\mathfrak{F}}=1$ or $N^{\mathfrak{F}}=N$. The case $N^{\mathfrak{F}}=1$ is impossible,
 since  $N$ is non-abelian, and $\mathfrak{F}\subseteq\mathfrak{S}$.
 Therefore $N^{\mathfrak{F}}=N$. By [30, proposition~A.4.13(a)] 
 $N$ is a direct product of subgroups, each isomorphic with a fixed simple non-abelian group.
 If  $N_{G}(G_{p})\cap N=N$, then $N_{p}=G_{p}\cap N\unlhd N$. By [30, proposition~A.4.13(b)] 
 $N_{p}$ is the direct product of a subset of the non-abelian factors of $N$. This is the contradiction.
  If $N_{G}(G_{p})\cap N\not=N$, then there is maximal subgroup $M$ of $N$ such that $N_{G}(G_{p})\cap N\leq M$ and $N^{\mathfrak{F}}\leq M$. We have the contradiction $N=N^{\mathfrak{F}}\leq M\not=N$.

\medskip

2. $N$  is an abelian $p$-group, $p$ is some prime. From  $G/N\in\mathfrak{F}\subseteq\mathfrak{S}$ and $N\in\mathfrak{S}$ it follows that $G$ is solvable.  From the uniqueness of $N$ and $\Phi(G)=1$ we conclude that $G=N\leftthreetimes M$, where $G^{\mathfrak{F}}=N=C_{G}(N)=F(G)$  and $M$ is a maximal subgroup of $G$, and moreover, $M\in\mathfrak{F}\subseteq \mathfrak{L}_{a}(1)$.

 Suppose that $M$  is nilpotent. By Lemma 1.4 $O_{p}(M)=1$, therefore $p\cap\pi(M)=\emptyset$. It follows that $M$ contains a normal Sylow $q$-subgroup  $M_{q}$ for some $q\in\pi(M)$ and $q\neq p$. Therefore $M_{q}= G_{q}$ is a Sylow $q$-subgroup of the group $G$.   From the uniqueness of $N$
it follows that $N_{G}(G_{q})\neq G$. Since $M$ is a maximal subgroup of $G$ and $M\subseteq
N_{G}(G_{q})$, we have $M=N_{G}(G_{q})$. But this contradicts the fact that $N_{G}(G_{q})$ $\mathfrak{F}$-$\mathrm{sn}$ $G$.

We assume that $M$ is non-nilpotent. Let  $\pi(G)=\{p_{1}, p_{2},\ldots, p_{n}\}$, where  $p_{1}=p$.  Consider the following cases.

\medskip

i) Let  $n=2$. Then $p\in\pi(M)$. By Lemma 1.4  $O_{p}(M)=1$. Since  $M\in\mathfrak{L}_{a}(1)$, by Lemma~3.2  $M\in\mathfrak{N}^{2}$. Therefore $M/F(M)$
 is nilpotent. We note that  $F(M)$  is a $p_{2}$-group. If $Q\in {\rm Syl}_{p_{2}}(M)$, then $Q$ is a normal subgroup of $M$, moreover, $Q\in {\rm Syl}_{p_{2}}(G)$ and $N_{G}(Q)=M$. By hypothesis
 $N_{G}(Q)=M$ $\mathfrak{F}$-$\mathrm{sn}$ $G$. Therefore $N=G^{\mathfrak{F}}\subseteq M$ and $G=NM\subseteq M$. This is the contradiction.

\medskip

ii) Let  $n\geq 3$.

We will to show that $N$ is a Sylow $p$-subgroup of $G$.
By Hall's theorem $G=G_{1}G_{2}\cdots G_{n}$, where  $G_{1}$, $G_{2}, \ldots ,  G_{n}$ are pairwise permutable Sylow  $p_{1}$-, $p_{2}$-, $\ldots$ , $p_{n}$-subgroups of $G$, respectively.
Let  $A_{i} = G_{1}G_{i}$, where  $i\not=1$. Since
 $|A_{i}|<|G|$, $N_{G}(G_{1})\cap A_{i}=N_{A_{i}}(G_{1})$ $\mathfrak{F}$-$\mathrm{sn}$ $A_{i}$ and $N_{G}(G_{i})\cap A_{i}=N_{A_{i}}(G_{i})$ $\mathfrak{F}$-$\mathrm{sn}$ $A_{i}$, we have $A_{i}\in\mathfrak{F}$.
 From $|\pi(A_{i})|=2$ by Lemma 3.2 it follows that  $A_{i}\in\mathfrak{N}^{2}$.
We note that $N\subseteq A_{i}$. Since  $N=C_{G}(N)$ and $p_{1}=p$, we see that $F(A_{i})$ is a $p$-group. From $A_{i}\in\mathfrak{N}^{2}$ it follows that $A_{i}/F(A_{i})\in\mathfrak{N}$. Then  $G_{1}/F(A_{i})\unlhd A_{i}/F(A_{i})$ and $G_{1}\unlhd A_{i}$. So, $G_{i}\subseteq A_{i}\subseteq N_{G}(G_{1})$. Hence  $G\subseteq N_{G}(G_{1})$ and $G_{1}\unlhd G$.
From $G_{1}\cap M\unlhd M$ and $O_{p}(M)=1$ it follows that $G_{1}\cap M=1$. So $G_{1}=N\in{\rm Syl}_{p}(G)$.

Thus $M$ is a $p'$-Hall subgroup of $G$. Let  $i\in\{2,\ldots, n\}$ and $S\in {\rm Syl}_{p_{i}}(M)$. Then  $S\in {\rm Syl}_{p_{i}}(G)$ and $N_{G}(S)\neq M$. We note that $N_{G}(S)\neq G$ because  $N=C_{G}(N)$ and $N$ is a $p$-group, $p\neq p_{i}$.

We will to show that  $N_{G}(S)\in\mathfrak{F}$.

Suppose that $N_{G}(S)\cap N=1$. Since $G/N\in\mathfrak{F}$ and $\mathfrak{F}$ is a hereditary formation, it follows that
 $N_{G}(S)N/N\cong N_{G}(S)/N_{G}(S)\cap N\cong N_{G}(S)\in\mathfrak{F}$.

Suppose now that  $N_{G}(S)\cap N=D\neq 1$. Then  $D\unlhd N_{G}(S)$ and $S\unlhd N_{G}(S)$. We have  $S\times D\unlhd N_{G}(S)$ and
 $N_{G}(S)=(S\times D)\leftthreetimes R$, where  $R$ is a $\{p_{1}, p_{i}\}'$-Hall subgroup of $N_{G}(S)$. From $G\in\mathfrak{S}$ by Hall’s theorem we deduce that $SR\leq M^{x}$ for some $x\in G$ and there is a $\{p_{i}\}'$-Hall subgroup $H$ from $G$ such that $DR\leq H$. From  ${\rm Syl}(H)\subseteq{\rm Syl}(G)$ it follows that $N_{G}(L)$ $\mathfrak{F}$-$\mathrm{sn}$ $G$ for any $L\in{\rm Syl}(H)$. By Lemma~1.7(1) $N_{H}(L)=N_{G}(L)\cap H$ $\mathfrak{F}$-$\mathrm{sn}$ $H$. Then  $H\in\mathfrak{F}$ by the choice of $G$. We note that $M^{x}\cong M\in\mathfrak{F}$. Since $\mathfrak{F}$ is hereditary we have  $N_{G}(S)/D\cong SR\in\mathfrak{F}$ and $N_{G}(S)/S\cong DR\in\mathfrak{F}$. We obtain $N_{G}(S)/S\cap D\cong N_{G}(S)\in\mathfrak{F}$.

 Consider  $T=NN_{G}(S)$.
From Lemma~1.7(1) $N_{G}(S)$ $\mathfrak{F}$-$\mathrm{sn}$ $T$.
 By theorem 15.10 [19] 
 $T\in\mathfrak{F}$. Let  $h$ be the maximal inner local screen formation $\mathfrak{F}$.
 By Lemma~1.5 [19] 
 it follows that $T/F_{p}(T)\in h(p)$. Because $N\leq F_{p}(T)$ and $N=C_{G}(N)$, we have
 $O_{p'}(T)=1$ and $N=F_{p}(T)$. Therefore $T/N\in h(p)$.
  Then $N_{G}(S)N/N\cong N_{G}(S)/N_{G}(S)\cap N\in h(p)$. Since  $\mathfrak{F}$ is
  a hereditary formation, it follows that  $h(p)$ is a hereditary
  formation, by the theorem~4.7 [19]. 
 Then    $(N_{G}(S)\cap M)N/N\cong N_{G}(S)\cap M/N_{G}(S)\cap N\cap M\cong N_{G}(S)\cap M\in h(p)$. We note that $N_{G}(S)\cap M=N_{M}(S)$. Therefore $N_{M}(S)\in h(p)$. By Lemma 3.4  $M\in h(p)$.
 Then  $G/F_{p}(G)\cong M\in h(p)$. By Lemma 1.5 
 $G\in\mathfrak{F}$, which contradicts the choice of $G$.
$\Box$

\medskip

{\bf Corollary~3.4.1} [13].
 {\it If the normalizers of all Sylow subgroups of a group $G$ are $\mathbb{P}$-subnormal, then $G$ is supersoluble.}

\medskip

{\bf Corollary~3.4.2} [29].
 {\it A group $G\in\mathfrak{N}^{2}$ if and only if all its Sylow subgroups are strongly ${\rm K}$-$\mathfrak{N}^{2}$-subnormal in $G$.}

\medskip

 {\bf Corollary~3.4.3} [29].
 {\it A group $G\in\mathfrak{NA}$ if and only if all its Sylow subgroups are strongly ${\rm K}$-$\mathfrak{NA}$-subnormal in $G$.}

\medskip

 {\bf Corollary~3.4.4.}
 {\it A group $G\in\mathfrak{L}_{a}(1)$ if and only if all its Sylow subgroups are strongly ${\rm K}$-$\mathfrak{L}_{a}(1)$-subnormal in $G$.}

\medskip

{\bf Remark~3.5.} Note that ${\rm w}^{*}_{\pi}\mathfrak{F}\subseteq\overline{{\rm W}}_{\pi}\mathfrak{F}$.
From [23, 25] it follows that $\overline{{\rm W}}\mathfrak{N}^{2}={\rm w}\mathfrak{N}^{2}=\mathfrak{S}$. But ${\rm w}^{*}\mathfrak{N}^{2}=\mathfrak{N}^{2}$.

\newpage

{\parindent=0mm
\textbf{\bf {\large References}}
}

\medskip

[1]
M.\,G. Bianci, A. Gillio Berta Mayri and P. Hauck.
On finite soluble groups with nilpotent Sylow normalizers.
{\it Arch. Math.} {\bf 47} (1986) 193--197.

[2]
V. Fedri and L. Serena.
Finite soluble groups with supersoluble Sylow normalizers.
{\it Arch. Math.}
{\bf 50} (1988) 11--18.

[3]
R.\,A. Bryce, V. Fedri and L. Serena.
Bounds on the Fitting length of finite soluble groups with supersoluble Sylow normalizers.
{\it Bull. Austral. Math. Soc.}
{\bf 44} (1991) 19--31.

[4]
A. Ballester-Bolinches and L.\,A. Shemetkov.
On normalizers of Sylow subgroups in finite groups.
{\it Siberian Math. J.} {\bf 40}(1) (1999) 1--2.

[5]
A. D'Aniello, C. De Vivo and G. Giordano.
Finite groups with primitive Sylow normalizers.
{\it Bolletino U.M.I.}
{\bf 8}(5-B) (2002) 235--245.

[6]
A. D'Aniello, C. De Vivo and G. Giordano.
Saturated formations and Sylow normalizers.
{\it Bull. Austral. Math. Soc.}
{\bf 69} (2004) 25--33.

[7]
A. D'Aniello, C. De Vivo, G. Giordano and M.\,D. P\'{e}rez-Ramos.
Saturated formations closed under Sylow normalizers.
{\it Bull. Austral. Math. Soc.}
{\bf 33} (2005) 2801--2808.

[8]
L. Kazarin, A. Mart\'{\i}nez-Pastor and M.\,D. P\'{e}rez-Ramos.
On the Sylow graph of a group and Sylow normalizers.
{\it Israel J. Math.}
{\bf 186} (2011) 251--271.

[9]
L. Kazarin, A. Mart\'{\i}nez-Pastor and M.\,D. P\'{e}rez-Ramos.
On Sylow normalizers of finite groups.
{\it J. Algebra Appl.}
{\bf 13}(3) (2014) 1350116--1--20.

[10]
G. Glaubermann.
Prime-power factor groups of finite groups II.
{\it Math. Z.}
{\bf 117} (1970) 46--56.

[11]
A.\,F. Vasil'ev, T.\,I. Vasil'eva and V.\,N. Tyutyanov.
On the finite groups of supersoluble type.
{\it Siberian Math. J.}
{\bf 51}(6) (2010) 1004--1012.

[12]
A.\,F. Vasil'ev, T.\,I. Vasil'eva and V.\,N. Tyutyanov.
On $\rm{K}$-$\mathbb{P}$-subnormal subgroups of finite groups.
{\it Math. Notes.} {\bf 95}(4) (2014) 471--480.

[13]
V.\,N. Kniahina, V.\,S. Monakhov.
On supersolvability of finite groups with
$\mathbb{P}$-subnormal subgroups.
{\it Internat. J. of Group Theory}  {\bf 2}(4) (2013) 21--29.

[14]
I. Zimmermann.
Submodular subgroups in finite groups.
{\it Math. Z.}
{\bf 202} (1989) 545--557.

[15].
R. Schmidt.
{\it Subgroup Lattices of Groups}
(Walter de Gruyter, 1994). 

[16]
V.\,A. Vasilyev.
Finite groups with submodular Sylow subgroups.
{\it Siberian Math. J.}
{\bf 56}(6) (2015) 1019--1027.

[17]
V.\,A. Vasilyev.
On the influence of submodular subgroups on the structure of finite groups.
{\it Vestnik Vitebsk Univ.}
{\bf 91}(2) (2016) 17--21. (In Russian)

[18]
T. Hawkes.
On formation subgroups of a finite soluble group.
{\it J. London Math. Soc.} {\bf 44} (1969) 243--250.

[19]
L.\,A. Shemetkov {\it Formations of finite groups} (Nauka, Moscow, 1987). (In Russian)

[20]
A. Ballester-Bolinches and L.\,M. Ezquerro.
{\it Classes of Finite Groups} (Springer, 2006).  

[21]
O.\,H. Kegel.
Untergruppenverb\"ande endlicher Gruppen,
die den Subnormalteilerverband echt enthalten.
{\it Arch. Math.}
{\bf 30}(3) (1978) 225--228.

[22]
A.\,F. Vasil'ev.
On the influence of primary $\mathfrak{F}$-subnormal
subgroups on the structure of the group.
{\it Voprosy Algebry (Problems in Algebra)}.
{\bf 8} (1995) 31--39. (In Russian)

[23]
A.\,F. Vasil'ev, T.\,I. Vasil'eva and A.\,S. Vegera.
Finite groups with generalized subnormal embedding of Sylow subgroups.
{\it Siberian Math. J.}
{\bf 57}(2) (2016) 200--212.

[24]
T.\,I. Vasil'eva and A.\,I. Prokopenko.
Finite groups with generally subnormal subgroups.
{\it Proceedings of the National Academy of Sciences of Belarus.} Series of Physical-Mathematical Sciences.
{\bf 3} (2006) 25--30. (In Russian)

[25]
A.\,F. Vasil'ev and T.\,I. Vasil'eva.
On finite groups with generally subnormal Sylow subgroups.
{\it PFMT} {\bf 4}(9) (2011) 86--91. (In Russian)

[26]
A.\,S. Vegera.
On local properties of the formations of groups with $\rm{K}$-$\mathfrak{F}$-subnormal Sylow subgroups.
{\it PFMT} {\bf 3}(20) (2014) 53--57. (In Russian)

[27]
V.\,S. Monakhov and I.\,L. Sokhor.
Finite groups with formation subnormal primary subgroups.
{\it Siberian Math. J.} {\bf 58}(4) (2017) 851--863.

[28]
V.\,I. Murashka.
Finite groups with given sets of $\mathfrak{F}$-subnormal subgroups.
{\it Asian-European J. Math.}
(2019) 2050073 (13 pages).  DOI: 10.1142/S1793557120500734.

[29]
A.\,F. Vasil'ev.
Finite groups with strongly $\rm{K}$-$\mathfrak{F}$-subnormal Sylow subgroups.
{\it PFMT}{\bf 4}(37) (2018) 66--71. (In Russian)

[30]
K. Doerk and T. Hawkes.
{\it Finite soluble groups}
(Walter de Gruyter, Berlin, New York, 1992). 

[31]
V.\,N. Semenchuk.
Minimal non $\mathfrak{F}$-subgroups
{\it Algebra and Logik} {\bf 18}(3) (1979) 348--382.

\bigskip

A.\,F.\,Vasil'ev, A.\,G.\,Melchenko

Francisk Skorina Gomel State University,
Gomel, Belarus.

E-mail address: formation56@mail.ru, melchenkonastya@mail.ru

\medskip

Т.\,I.\,Vasil'eva

Francisk Skorina Gomel State University,
Belarusian State University of Transport,
Gomel, Belarus.

E-mail address: tivasilyeva@mail.ru





\end{document}